\documentclass[11pt,a4,fleqn]{article}
\usepackage{graphicx}
\usepackage{amsmath,amssymb,latexsym,graphics,epsfig}
\usepackage{hyperref}
\usepackage{color}
\usepackage{amsthm}

\newtheorem{theorem}{\bf Theorem}[section]

\newtheorem{lemma}[theorem]{\bf Lemma}

\newtheorem{emp}[theorem]{\bf Claim }
\newtheorem{remark}[theorem]{\bf Remark}

\newtheorem{conjecture}[theorem]{\bf Conjecture}
\newtheorem{fact}[theorem]{\bf Fact}

\newcommand{\epf}{\hfill\hbox{\rule{3pt}{6pt}}\\}

\numberwithin{equation}{section}

\begin{document}
\title{{\Large A proof for a conjecture of Gy\'{a}rf\'{a}s, Lehel, S\'{a}rk\"{o}zy and Schelp on Berge-cycles}}
\author{G.R. Omidi  \\[2pt]
{\small  Department of Mathematical Sciences, Isfahan University of Technology},\\
{\small Isfahan, 84156-83111, Iran}\\
{\small School of Mathematics, Institute for Research in Fundamental Sciences (IPM),}\\
{\small P.O. Box 19395-5746, Tehran, Iran }\\[2pt]
{romidi@cc.iut.ac.ir}}

\date{}

\maketitle \footnotetext[1] {This research is partially
carried out in the IPM-Isfahan Branch and in part supported
by a grant from IPM (No. 92050217).} \vspace*{-0.5cm}

\begin{abstract} \rm{} It has been conjectured that for any fixed $r\geq 2$ and sufficiently large $n$, there is a monochromatic Hamiltonian Berge-cycle in every $(r-1)$-coloring of the edges of $K_{n}^{r}$, the complete $r$-uniform hypergraph on $n$ vertices. In this paper we prove this conjecture.

\noindent{\small { Keywords:} Monochromatic Hamiltonian  Berge-cycle, Colored complete uniform hypergraphs.}\\
{\small AMS subject classification: 05C65, 05C45, 05D10.}

\end{abstract}

\section{\normalsize Introduction}

For a given $r \geq 2$ {and $n\geq r$}, an {\it $r$-uniform Berge-cycle of length $n$}, denoted by $C_n^{r}$, is an $r$-uniform hypergraph with the core sequence $v_{1}, v_{2}, \ldots ,v_{n}$ as the vertices, and distinct edges $e_{1},e_{2},\ldots,e_{n}$ such that $e_{i}$ contains $v_{i}, v_{i+1}$, where addition in indices is in modulo $n$. {The case $r=2$ gives the usual definition of the cycle $C_n$ on $n$ vertices in graph case.} A Berge-cycle of length $n$ in a hypergraph with $n$ vertices is called a {\it Hamiltonian Berge-cycle}.\\

For an $r$-uniform hypergraph $H$, the {\it Ramsey number} $R_{k}(H)$ is the minimum integer $n$ such that there is a monochromatic copy of $H$ in every $k$-edge coloring of $K_{n}^{r}$. The existence of such a positive integer is guaranteed by the Ramsey's classical result in \cite{Ramsey}. Recently, the Ramsey numbers of various variations of cycles in uniform hypergraphs have been studied, e.g. see \cite{hax1, hax2, OS}. In this regard, Gy\'{a}rf\'{a}s et al. proposed the following conjecture for Berge-cycles:

\begin{conjecture}{\rm \cite{gyarf1}}\label{GLSS}
 Assume that $r\geq 2$ is fixed and $n$ is sufficiently large. Then every $(r-1)$-edge coloring of $K_{n}^{r}$ contains a monochromatic Hamiltonian Berge-cycle.
\end{conjecture}

Conjecture \ref{GLSS} states that for a given $r\geq 2$ we have $R_{r-1}(C_n^{r})=n$ when $n$ is sufficiently large. The case $r=2$ is trivial {since every complete graph $K_n$ has a Hamiltonian cycle.} The case $r=3$ was proved {by Gy\'{a}rf\'{a}s et al.} in \cite{gyarf1}. Recently,  Maherani and the author gave a proof for the case $r=4$; see \cite{MO}. For general $r$, the asymptotic form of Conjecture \ref{GLSS} was proved {by A. Gy\'{a}rf\'{a}s, G.N. S\'{a}rk\"{o}zy and E. Szemer\'{e}di  using the method of Regularity Lemma; see \cite{GSS}.} To see more results on Conjecture \ref{GLSS}, we refer the reader to \cite{gyarf1,gyarf2,GSS} and references therein. In this paper, we establish  Conjecture \ref{GLSS}. Based on the above results on this conjecture it only suffices to give a proof for $r\geq 5$. The main result of this paper is the following theorem:
\begin{theorem}\label{main result}
Suppose that $r\geq 4$ and $n>6r{4r \choose r-1}$. Then in every $(r-1)$-edge coloring of $K_{n}^{r}$ there is a monochromatic Hamiltonian Berge-cycle.
\end{theorem}

For a given $r\geq 2$, let $p(r)$ be the minimum value of $m$ for which the statement of Conjecture \ref{GLSS} holds for any $n\geq m$.
Theorem \ref{main result} guarantees the existence of such a function $p(r)$ (in fact it shows that $p(r)\leq 6r{4r \choose r-1}+1$). Determining $p(r)$ seems to be an interesting problem, though we will not make any serious attempt in this direction. At present, we do not know much about $p(r)$. Our conjecture is that $p(r)$ is much less than $6r{4r \choose r-1}+1$, at least for small values of $r$. An indication for this is given by $p(3)=5$ (see \cite{gyarf1}) and $p(4)\leq 85$ (see \cite{MO}). {In the rest of this paper, for a real number $x$ by $\lfloor x\rfloor$ (resp. $\lceil x\rceil$) we mean the greatest integer not exceeding $x$ (resp. the least integer not less than $x$).}

\section{\normalsize  Basic definitions and some preliminaries}

Before we give our proof we present some definitions. Assume that $H$ is an $r$-uniform hypergraph. The {\it shadow graph} $\Gamma (H)$ is a graph with vertex set $V(H)$, where two vertices are adjacent if they are covered by at least one edge of $H$. Consider an $(r-1)$-edge coloring of $H=K_{n}^{r}$ with colors $1,2,\ldots,r-1$ and assume that $G=\Gamma (H)$ (so $G$ is a complete graph). For each edge $e=xy$ of $G$, we assign a list $L(e)$ of colors of all edges of $H$ containing $x$ and $y$. For an edge $e\in E(G)$, the color $i\in L(e)$ is {\it good}, if at least $r-1$ edges (of $H$) of color $i$ contain both vertices of $e$.
We consider a new multi-coloring $L^{*}$ for the edges of $G$. For each edge $e\in E(G)$, assume that $L^{*}(e)\subseteq L(e)$ is the set of all good colors for $e$. Throughout this paper, for each natural number $m$, assume that $[m]=\{1,2,\ldots,m\}$. For each vertex $x\in V(G)$ and any $1\leq i\leq r-1$ assume that $$U_i(x)=\{y\in V(G)\setminus \{x\}|i\in L^{*}(xy)\}, \overline{U}_i(x)=\{y\in V(G)\setminus \{x\}|i\notin L^{*}(xy)\},$$ and $d_i(x)$ is the number of edges of color $i$ containing $x$ in $H$. For any $I\subseteq [r-1]$, set $U_I(x)=\bigcap_{i\in I}U_i(x)$ and $\overline{U}_I(x)=\bigcap_{i\in I}\overline{U}_i(x)$. We say that a set of vertices $S\subseteq V(G)$ {\it avoids} the set of colors $W\subseteq [r-1]$, if for each $i\in W$ there is a vertex $x\in S$ with $d_i(x)\leq {4r \choose r-1}$ or an edge $e=xy$ for $x,y\in S$ with $i\notin L^{*}(e)$. We will use the following lemmas in the proof of Theorem \ref{main result}.

{
\begin{lemma}\label{HB0}\cite{gyarf2}
 Assume that $r\geq3$ and $H=K_{n}^{r}$ is an $(r-1)$-edge colored complete $r$-uniform hypergraph on $n$ vertices. Also suppose that $G=\Gamma (H)$ and there is a monochromatic Hamiltonian cycle in $G$ under multi-coloring $L^{*}$. Then there is a monochromatic Hamiltonian  Berge-cycle in $H$.
\end{lemma}
}

\begin{lemma}\label{H1}\cite{Bondy}
Let $G$ be a simple graph and let $u$ and $v$ be nonadjacent vertices in $G$ such that $d_G(u)+d_G(v)\geq n.$ Then $G$ is Hamiltonian if and only if $G+uv$ is Hamiltonian.
\end{lemma}

\begin{lemma}\label{Ch1}\cite{Bondy}
Let $G$ be a simple graph with degree sequence $0\leq d_1\leq d_2\leq \cdots\leq d_n<n$ and $n\geq 3$. If for each $i<n/2$, we have $d_i>i$ or $d_{n-i}\geq n-i$, then $G$ is Hamiltonian.
\end{lemma}

{The following simple remark can be proved by induction on $m$ and it will be used later on.}
{\begin{remark}\label{remark1}
Assume that $a_m\geq a_{m-1}\geq\cdots \geq a_1\geq a>0$ are real numbers and $a_1+\cdots +a_m=l$. Then $\prod_{i=1}^{m}a_i\geq a^{m-1}(l-(m-1)a)$.
\end{remark}}

\section{\normalsize  Outline of the proof}

Here, we sketch the main ideas of our proof for Theorem \ref{main result}. Suppose to the contrary that there is no monochromatic Hamiltonian Berge-cycle in a given $(r-1)$-edge coloring $c$ of $H=K_{n}^{r}$ with colors $1,2,\ldots,r-1$.
{We will show that (see Claim \ref{claime1})}, by suitably renaming of colors, for some $0\leq f\leq r-2$ there are distinct vertices $x$ and $\{y_i\}_{i=1}^{r-1}$ such that $|\overline{U}_{r-1}(x)|\geq (n-1)/2$, $i\notin L^{*}(xy_i)$ for any $f+1\leq i\leq r-1$ and $\{y_i\}_{i=1}^{f}$ avoids $[f]$. We choose distinct vertices $x$ and $\{y_i\}_{i=1}^{r-1}$ with these properties and maximum $f$. Without loss of generality we assume that,
$$|\overline{U}_{f+1}(x)|\leq |\overline{U}_{f+2}(x)|\leq \cdots \leq |\overline{U}_{r-1}(x)|.$$ Then we divide our proof to some cases and in each case using the distinct vertices $x$ and $\{y_i\}_{i=1}^{r-1}$ we construct a new graph $\Gamma$ on $V(H)$ so that any Hamiltonian cycle in $\Gamma$ can be extended to a monochromatic Hamiltonian Berge-cycle of color $f+1$ in $H$. In fact $V(\Gamma)=V(H)$ and for any two adjacent vertices $u$ and $v$ of $\Gamma$, there exists an edge $f_{uv} \in E(H)$ of color $f+1$ containing $u$ and $v$. Moreover $f_{uv}\neq f_{u'v'}$ for almost any two distinct edges $uv$ and $u'v'$ in $E(\Gamma)$. {In overall $\Gamma$ can be defined as follows. The choices of the vertices  $D=\{y_i\}_{i=1}^{r-1}\cup \{x\}$ imply that for almost all vertices $u\in V(H)$ there are so many vertices $v$ such that there is an edge $e_{uv}$ in $H$ of color $f+1$ containing $u,v$ with $|e_{uv}\cap D|\geq r-2$. Now we consider the new graph with vertex set $V(H)$ and edges $uv$ mentioned above. Then we add a few suitable edges (the edges $E_3$ and $E_4$ in Page 9) to this graph to get a new graph $\Gamma$ with minimum degree at least $2r+1$.} To complete our proof {(in fact to get a contradiction to our incorrect assumption) it suffices} to show that  $\Gamma$ is a Hamiltonian graph. To do this, we show that the degree sequence of the graph $\Gamma$ satisfies {the  Chv\'{a}tal's condition in Lemma \ref{Ch1}}. More precisely if $d_1\leq d_2\leq \cdots\leq d_n$ are degrees of the vertices of $\Gamma$, then for each $i\leq \frac{n}{2}$, we have $d_{i}>i$ or  $d_{n-i}\geq n-i$.  Hence by Lemma \ref{Ch1}, $\Gamma$ is Hamiltonian and we are done.

\section{The proof}

{\it Proof of Theorem \ref{main result}:} Suppose to the contrary that there is no monochromatic Hamiltonian Berge-cycle in a given $(r-1)$-edge coloring $c$ of $H=K_{n}^{r}$ with colors $1,2,\ldots,r-1$. For each $1\leq i\leq r-1$, let $W_i$ be the set of all edges $e$ of $G=\Gamma (H)$ for which $i\notin L^{*}(e)$. {Using Lemma \ref{HB0},} we may assume that {the subgraph of $G$ with vertex set $V(G)$ and edge set $E(G)\setminus W_i$} is not Hamiltonian. Now consider $S_i\subseteq W_i$ with minimum cardinality, such that the spanning subgraph of $G$ induced by {$E(G)\setminus S_i$} is not Hamiltonian. Assume that $G_i$ and $G^c_i$ are the spanning subgraphs of $G$ induced by $S_i$ and {$E(G)\setminus S_i$,} respectively. For each color $1\leq i\leq r-1$, respectively assume that $T_i$ and  $R_i$ are the sets of all isolated vertices and all vertices with degree at least $(n-1)/2$ of $G_i$. Also, assume that $Q_i=V(G_i)\setminus (T_i\cup R_i)$. {We need the following fact frequently in our proof.
\begin{fact}\label{Gi}
For each $1\leq i\leq r-1$, $G^c_i$ is non-Hamiltonian. Moreover for each $e\in E(G_i)$, we have $i\notin L^{*}(e)$ and $G^c_i+e$ is Hamiltonian.
\end{fact}}

For any two non-adjacent vertices $x$ and $y$ of $G^c_i$, by Fact \ref{Gi} the graph $G^c_i+xy$ is Hamiltonian and so, by Lemma \ref{H1}, we have $d_{G^c_i}(x)+d_{G^c_i}(y)\leq n-1$. {Therefore, we have the following fact on the sums of degrees of adjacent vertices in $G_i$.
\begin{fact}\label{Gi1}
For any two adjacent vertices $x$ and $y$ of $G_i$ we have  $d_{G_i}(x)+d_{G_i}(y)\geq n-1$.
\end{fact}
This fact implies that $Q_i$ is an independent set in $G_i$. If $R_i=\emptyset$ for some $i$, then since $Q_i$ is an independent set, the graph $G_i$ has no edge and so $G^c_i$ is a complete graph, a contradiction to the fact that $G^c_i$ is non-Hamiltonian. Hence $R_i\neq\emptyset$} (see Section 2 for the notations that are not defined here). Now we claim that $|R_i|\geq|T_i|$ for each $1\leq i\leq r-1$. Assume to the contrary that for some $i$ we have $|R_i|<|T_i|$. Let $R_i=\{x_1,x_2,\ldots,x_{|R_i|}\}$, $T_i=\{y_1,y_2,\ldots,y_{|T_i|}\}$ and $Q_i=\{z_1,z_2,\ldots,z_{|Q_i|}\}$. Obviously
$$C=y_1x_{1}\ldots y_{|R_i|}x_{|R_i|}y_{|R_i|+1}\ldots y_{|T_i|}z_1\ldots z_{|Q_i|},$$
is a Hamiltonian cycle in $G^c_i$, a contradiction. By the same argument, we have $|R_i\cup Q_i|>|T_i|$. {Therefore, we have the following fact.
\begin{fact}\label{Gi2}
For each $1\leq i\leq r-1$, we have
\begin{itemize}
\item[$\bullet$]
$Q_i$ is an independent set in $G_i$,
\item[$\bullet$]
$R_i\neq\emptyset$ and $|R_i|\geq|T_i|$,
\item[$\bullet$]
$|R_i\cup Q_i|>|T_i|$.
\end{itemize}
\end{fact}
}
An argument similar to the proof of Claim 2.3 of Theorem 2.2 in \cite{MO} (set $t=2$ and follow the proof) yields the following result:

\begin{emp}\label{losses}
Let $P\subseteq [r-1]$ and $|P|=p$. Then there is a set of vertices $Q\subseteq V(G)$ with $|Q|\leq p+1$ such that $Q$ avoids $P$.
\end{emp}

First assume that there is a subset $S\subseteq V(G)$ that avoids a set of colors containing at least $|S|+1$ colors $c_1,c_2,\ldots,c_{|S|+1}$. Using Claim \ref{losses}, there is a subset $S'\subseteq V(G)$ containing at most $r-1-|S|$ vertices that avoids $[r-1]\setminus \{c_1,c_2,\ldots,c_{|S|+1}\}$. Now $S\cup S'$ avoids $[r-1]$, which is impossible since the number of edges in $H$ containing $S\cup S'$ is ${n-|S\cup S'| \choose r-|S\cup S'|}\geq n-r+1$ and for each $1\leq i\leq r-1$ the number of edges of color $i$ containing $S\cup S'$ is at most ${4r \choose r-1}$ (note that $n>6r{4r \choose r-1}$).
Therefore, each subset $S\subseteq V(G)$ avoids at most $|S|$ colors in $[r-1]$.

\begin{emp}\label{claime1}
{By suitably renaming the colors, there are distinct vertices $x$ and $\{y_i\}_{i=1}^{r-1}$ such that $|\overline{U}_{r-1}(x)|\geq (n-1)/2$ and for some $0\leq f\leq r-2$, $\{y_i\}_{i=1}^{f}\subseteq \bigcap_{i=f+1}^{r-1}T_i$, the set of vertices $\{y_i\}_{i=1}^{f}$ avoids $[f]$ and $i\notin L^{*}(xy_i)$ for any $f+1\leq i\leq r-1$.}
\end{emp}

 {\it Proof of Claim \ref{claime1}:} Let $S=\{y_i\}_{i=1}^{f}\subseteq V(G)$ be the largest subset of vertices with $f\leq r-1$ that avoids a set containing $f$ colors. Note that it is possible to have $S=\emptyset$. Without any loss of generality, we may assume that $S$ avoids $[f]$. The case $f=r-1$ is impossible, since the number of edges in $H$ containing $S$ is $n-r+1> 6r{4r \choose r-1}-r+1$ and for each $1\leq i\leq r-1$ the number of edges of color $i$ containing $S$ is at most ${4r \choose r-1}$. Hence $f\leq r-2$.  If $y_i\notin T_j$ for some $1\leq i\leq f$ and $f+1\leq j\leq r-1$, then there is a vertex $v\in V(G)$ such that $j\notin L^{*}(vy_i)$ and so $S\cup\{v\}$ avoids $[f]\cup\{j\}$, a contradiction to the maximality of $S$. Hence
 $$~~~~~~~~~~~~~~~~~~~~~~~~~~~~~~~~S\subseteq \bigcap_{i=f+1}^{r-1}T_i.~~~~~~~~~~~~~~~~~~~~~~~~~~~~~(1)$$
 If $f=r-2$, then {choose} $x\in R_{r-1}$ and $y_{r-1}\in N_{G_{r-1}}(x)$. Since $d_{G_{r-1}}(x)\geq (n-1)/2$, we have $|\overline{U}_{r-1}(x)|\geq (n-1)/2$ and so there is nothing to prove. Now let $f\leq r-3$. {If for some $x\in V(G)$ and for some $f+1\leq i,j\leq r-1$ with $i\neq j$ we have $\overline{U}_i(x)\cap \overline{U}_j(x)\neq\emptyset$, then for any $v\in\overline{U}_i(x)\cap \overline{U}_j(x)$ the set $S\cup\{x,v\}$ avoids $[f]\cup\{i,j\}$, a contradiction to the maximality of $f$. Hence the following fact holds.
 \begin{fact}\label{Gi3}
For any $f+1\leq i,j\leq r-1$ with $i\neq j$ and for each $x\in V(G)$, we have
 $\overline{U}_i(x)\cap \overline{U}_j(x)=\emptyset$
\end{fact}}

 {Now we claim that there is a vertex  $x\in\bigcup_{i=f+1}^{r-1}R_i\setminus\bigcup_{i=f+1}^{r-1}T_i$. If there is such a vertex $x$, then the proof of Claim \ref{claime1} will be finished by an easy argument. To see this, without any loss of generality assume that $x\in R_{r-1}$. Since $x$ has degree at least $(n-1)/2$ in $G_{r-1}$, we have $|\overline{U}_{r-1}(x)|\geq (n-1)/2$. On the other hand, for each $i=f+1,\ldots, r-1$, we have  $x\in R_i\cup Q_i$. Hence for each $f+1\leq i\leq r-1$, there is a vertex $y_i$ with $xy_i\in E(G_i)$ {and so using Fact \ref{Gi} we have} $i\notin L^{*}(xy_i)$. Therefore, the vertices $x$ and $\{y_i\}_{i=1}^{f}$ have the desired properties in Claim \ref{claime1} and we are done. Now to show that $\bigcup_{i=f+1}^{r-1}R_i\setminus\bigcup_{i=f+1}^{r-1}T_i\neq\emptyset$ with a contrary assume $$~~~~~~~~~~~~~~~~~~~~~~~~~~~~~~~~\bigcup_{i=f+1}^{r-1}R_i\subseteq\bigcup_{i=f+1}^{r-1}T_i.~~~~~~~~~~~~~~~~~~~~~~~~~~~~~(2)$$}
We consider the following cases and in each case we get a contradiction.

\bigskip\noindent \textbf{Case 1.} $R_i\cap R_j=\emptyset$ for any $f+1\leq i,j\leq r-1$.\\
{By Fact \ref{Gi2} for each $i\leq r-1$, we have $|R_i|\geq|T_i|$. On the other hand, we have $R_i\cap R_j=\emptyset$ for any $f+1\leq i,j\leq r-1$ and using (2), $\bigcup_{i=f+1}^{r-1}R_i\subseteq\bigcup_{i=f+1}^{r-1}T_i$. Therefore, we have} $|R_i|=|T_i|$ for each $f+1\leq i\leq r-1$, $\bigcup_{i=f+1}^{r-1}R_i=\bigcup_{i=f+1}^{r-1}T_i$ and $T_i\cap T_j=\emptyset$ for any $f+1\leq i,j\leq r-1$ and $i\neq j$. {Note that by (1) we have $S\subseteq \bigcap_{i=f+1}^{r-1}T_i$, and therefore $f=0$. Using Fact \ref{Gi2} for each $1\leq i\leq r-1$, we have $R_i\neq\emptyset$. On the other hand, for each $1\leq i\leq r-1$ we have $|R_i|=|T_i|$ and the degree of each vertex of $R_i$ in $G_i$ is at least $(n-1)/2$.  Hence for each $1\leq i\leq r-1$, $Q_i\neq\emptyset$.} For each $1\leq i\leq r-1$, we have $d_{G_i}(w)\leq n-1-|T_i|$ when $w\in R_i$, and $d_{G_i}(w)\leq |R_i|$ when $w\in Q_i$. On the other hand, $|R_i|=|T_i|$ and {by Fact \ref{Gi1} we have $d_{G_i}(x)+d_{G_i}(y)\geq n-1$ for any two adjacent vertices $x$ and $y$ of $G_i$.} Therefore, for each $i$, the bipartite subgraph of $G_i$ with color classes $R_i$ and $Q_i$ is complete, and also the subgraph of $G_i$ induced by $R_i$ is a complete graph. Without any loss of generality, suppose that for every $1\leq i\leq r-2$, we have $|R_{r-1}|\leq|R_{i}|$. {Now for every $1\leq i\leq r-2$, set $A_i=R_{r-1}\cap T_{i}$ and $B_i=R_{r-1}\cap Q_{i}=R_{r-1}\setminus A_i$ (note that $R_{r-1}\cap R_i=\emptyset$).} Also, with no loss of generality, assume that  $|A_i|\leq|A_j|$ for {$i\leq j\leq r-1$.}

{First assume $A_{r-3}$ is non-empty. Clearly $R_{t}\setminus T_{r-1}$ is non-empty for some $t\in\{r-3,r-2\}$, since $|T_{r-1}|=|R_{r-1}|<|R_{r-2}\cup R_{r-3}|$. We claim that $B_{t}\neq\emptyset$. To see this, first suppose that $t=r-3$. If $B_{r-3}=\emptyset$, then $R_{r-1}=A_{r-3}\subseteq T_{r-3}$ and so
$R_{r-1}=A_{r-2}\subseteq T_{r-2}$ (note that $|A_{r-3}|\leq|A_{r-1}|$ and $A_{r-3}\cup A_{r-2}\subseteq R_{r-1}$). Hence $T_{r-3}\cap T_{r-2}\cap R_{r-1}=R_{r-1} \neq\emptyset$, a contradiction to the fact that $T_i\cap T_j=\emptyset$ for any $1\leq i,j\leq r-1$ and $i\neq j$. Now suppose that $t=r-2$. If $B_{r-2}=\emptyset$, then $R_{r-1}=A_{r-2}\subseteq T_{r-2}$ and so $A_{r-3}\subseteq T_{r-2}\cap T_{r-3}$, again a contradiction to the fact that $T_i\cap T_j=\emptyset$ for any $1\leq i,j\leq r-1$ and $i\neq j$. Now for that $t$, choose two vertices $u\in B_{t}$ and $v\in R_{t}\setminus T_{r-1}$. Since $uv$ is an edge of $G_t$, using Fact \ref{Gi} we have $t\notin L^{*}(uv)$. On the other hand, $v\in R_{t}\setminus T_{r-1}$ and $R_t\cap R_{r-1}=\emptyset$ and so $v\in Q_{r-1}$. Therefore, $uv$ is an edge of $G_{r-1}$, and again using Fact \ref{Gi} we have $r-1\notin L^{*}(uv)$ and so $\{u,v\}$ avoids $\{t,r-1\}$, which contradicts the fact that $f=0$.}

 Now assume that $A_{r-3}=\emptyset$. Then $A_{1}=\cdots=A_{r-3}=\emptyset$, and therefore $R_{r-1}\subseteq T_{r-2}$, since $\bigcup_{i=f+1}^{r-1}R_i=\bigcup_{i=f+1}^{r-1}T_i$. If $R_{i}\setminus T_{r-1}$ is non-empty for some $i\in\{1,\ldots,r-3\}$, then $i,r-1\notin L^{*}(uv)$ for all $u\in B_{i}$ and $v\in R_{i}\setminus T_{r-1}$ {(note that since $R_{r-1}\cap R_i=\emptyset$ and $A_i=\emptyset$, we have $R_{r-1}=B_i\neq\emptyset$)} and so $\{u,v\}$ avoids $\{i,r-1\}$, which is impossible. {Otherwise, $\bigcup_{i=1}^{r-3} R_i\subseteq T_{r-1}$. On the other hand, $|T_{r-1}|=|R_{r-1}|\leq |R_i|$ for every $1\leq i\leq n-1$. Hence $r=4$ and}  $R_{1}=T_{3}$. Since $\bigcup_{i=1}^{3}R_i=\bigcup_{i=1}^{3}T_i$ and $A_1=\emptyset$, we have $R_{3}=T_{2}$ and $R_{2}=T_{1}$ and hence  $R_{1}\subseteq Q_{2}$, $R_{2}\subseteq Q_{3}$ and $R_{3}\subseteq Q_{1}$. {Now since for each $1\leq i\leq 3$ the bipartite subgraph of $G_i$ with color classes $R_i$ and $Q_i$ is complete, for any three vertices $v_i\in R_i$, where $i=1,2,3$, we have $v_1v_3\in E(G_1)$, $v_1v_2\in E(G_2)$, and  $v_2v_3\in E(G_3)$  and so using Fact \ref{Gi} we have $1\notin L^{*}(v_{1}v_{3})$, $2\notin L^{*}(v_{1}v_{2})$ and $3\notin L^{*}(v_{2}v_{3})$.} Therefore, $\{v_{1},v_{2},v_{3}\}$ avoids $[3]=\{1,2,3\}$, which is again impossible.

\bigskip\noindent \textbf{Case 2.} $R_i\cap R_j\neq\emptyset$ for some $f+1\leq i,j\leq r-1$ and $i\neq j$.\\
Without any loss of generality, assume that $R_{r-2}\cap R_{r-1}\neq\emptyset$, and let $x\in R_{r-2}\cap R_{r-1}$.
Note that $f\leq r-3$ and {using Fact \ref{Gi3}, we have $\overline{U}_{r-2}(x)\cap \overline{U}_{r-1}(x)=\emptyset$.} Therefore, $d_{G_{r-2}}(x)=d_{G_{r-1}}(x)=(n-1)/2$, $N_{G_{r-2}}(x)\cap N_{G_{r-1}}(x)=\emptyset$ and $V(G)=\{x\}\cup N_{G_{r-2}}(x)\cup N_{G_{r-1}}(x)$. Hence $T_{r-2}\cap T_{r-1}=\emptyset$ and so $f=0$ ({note that by (1) we have $S=\{y_i\}_{i=1}^{f}\subseteq \bigcap_{i=f+1}^{r-1}T_i$}). We remind that $r\geq 4$. {With no loss of generality assume that }$|R_{r-3}\cap N_{G_{r-1}}(x)|\geq|R_{r-3}\cap N_{G_{r-2}}(x)|$.

First assume $x\in R_{r-3}$. Then for each $y\in N_{G_{r-3}}(x)$ the set $\{x,y\}$ clearly avoids a set containing $r-3$ and one of the colors $r-2$ or $r-1$, a contradiction to the fact that $f=0$. In fact $\{x,y\}$  avoids $\{ r-3, r-2\}$ if $y\in N_{G_{r-2}}(x)$ and $\{x,y\}$  avoids $\{ r-3, r-1\}$ if $y\in N_{G_{r-1}}(x)$.

Now assume $x\notin R_{r-3}$. Then since $R_{r-3}\neq\emptyset$, $V(G)=\{x\}\cup N_{G_{r-2}}(x)\cup N_{G_{r-1}}(x)$ and $|R_{r-3}\cap N_{G_{r-1}}(x)|\geq|R_{r-3}\cap N_{G_{r-2}}(x)|$ we have $R_{r-3}\cap N_{G_{r-1}}(x)\neq\emptyset$. Now consider $y\in R_{r-3}\cap N_{G_{r-1}}(x)$. If there is a vertex $z\in N_{G_{r-3}}(y)\cap N_{G_{r-2}}(x)$, then $\{x,y,z\}$ avoids $\{r-3,r-2,r-1\}$, contradicting again $f=0$. Therefore, $N_{G_{r-3}}(y)\subseteq N_{G_{r-1}}(x)\cup\{x\}$. {Since $y\in N_{G_{r-1}}(x)$, $d_{G_{r-2}}(x)=d_{G_{r-1}}(x)=(n-1)/2$ and} $d_{G_{r-3}}(y)\geq (n-1)/2$ we have $x\in N_{G_{r-3}}(y)$ and so $\{x,y\}$ avoids $\{r-3,r-1\}$, which is impossible.\epf

We choose distinct vertices $x$ and $\{y_i\}_{i=1}^{r-1}$ with the desired properties mentioned in Claim \ref{claime1} and maximum $f$. In the sequel, for simplicity we denote $U_I(x)$ and $\overline{U}_I(x)$  (for $I\subseteq[r-1]$) by $U_I$ and $\overline{U}_I$, respectively. {Also, for simplicity we denote $U_{i}(x)$ and $\overline{U}_{i}(x)$ ($1\leq i\leq r-1$) by $U_i$ and $\overline{U}_i$, respectively. Using Claim \ref{claime1} we have $|\overline{U}_{r-1}|\geq(n-1)/2$ and by Fact \ref{Gi3}, $\overline{U}_i\cap \overline{U}_j=\emptyset$ for any $f+1\leq i,j\leq r-1$ with $i\neq j$. Hence $|U_{r-1}|\geq |U_i|$ for each $1\leq i\leq r-1$ and without loss of generality we may assume that,
$$~~~~~~~~~~~~~~~~~~~~~~~~~~~~~~~~|\overline{U}_{f+1}|\leq |\overline{U}_{f+2}|\leq \cdots \leq |\overline{U}_{r-1}|.~~~~~~~~~~~~~~~~~~~~~~~~~~~~~(3)$$}
{Also, by Claim \ref{claime1} we have $f\leq r-2$.} Let $Y=\{y_1,y_2,\ldots,y_{r-1}\}\setminus\{y_{f+1}\}$ {and $Y_i=Y\setminus\{y_{i}\}$ for every $1\leq i\leq r-1$. We will use the following simple fact in our proof. It follows from the fact that $\{y_i\}_{i=1}^{f}$ avoids $[f]$ and for each $f+1\leq i\leq r-1$ we have $i\notin L^{*}(xy_i)$.
\begin{fact}\label{avoid1}
For every $1\leq i\leq r-1$ and $i\neq f+1$, the set of vertices $Y_{i}\cup \{x\}$ avoids the set of colors $[r-1]\setminus\{i,f+1\}$. Also, $Y\cup \{x\}$ avoids $[r-1]\setminus\{f+1\}$.
\end{fact}}

{Also we need the following fact in our proof later on.
\begin{fact}\label{ui}
For every $1\leq i\leq r-1$ and $i\neq f+1$, we have $\overline{U}_{i}\cap (Y_{i}\cup\{y_{f+1}\})=\emptyset$. Moreover $\overline{U}_{f+1}\cap Y_{f+1}=\emptyset$.
\end{fact}
The proof of Fact \ref{ui} is trivial. In fact if for $i\neq f+1$ we have $\overline{U}_{i}\cap (Y_{i}\cup\{y_{f+1}\})\neq\emptyset$, then the set of vertices $Y_{i}\cup\{x,y_{f+1}\}$ avoids all colors $[r-1]$. But this is impossible, since the number of edges in $H$ containing $Y_{i}\cup\{x,y_{f+1}\}$ is $n-r+1> 6r{4r \choose r-1}-r+1$ and for each $1\leq i\leq r-1$ the number of edges of color $i$ containing $Y_{i}\cup\{x,y_{f+1}\}$ is at most ${4r \choose r-1}$. The proof for the second result in  Fact \ref{ui} is similar.}

In the rest of our proof, we define a Hamiltonian graph $\Gamma$ with $V(\Gamma)=V(H)$, in such a way that every Hamiltonian cycle $C$ of $\Gamma$  can be extended to a monochromatic Hamiltonian Berge-cycle of $H$. For this, we consider the following cases:

\bigskip\noindent \textbf{Case 1.} $f=r-2.$\\
Consider a graph $\Gamma$ with vertex set $V(\Gamma)=V(H)$ and edge set
$E(\Gamma)=E_{1}\cup E_{2}$, where $E_i$'s are defined as follows:
\begin{itemize}
\item[$\bullet$]
$E_{1}=\{ uv\vert u,v\in V(\Gamma)\setminus Y, c(Y\cup\{u,v\})=r-1\}$. For each $uv\in E_{1}$, set $e_{uv}=Y\cup\{u,v\}$ and $F_{1}=\{e_{uv}|uv\in E_{1}\}$.
\item[$\bullet$]
$E_{2}=\{ y_iv\vert 1\leq i\leq r-2, v\in V(\Gamma)\setminus Y,\}$.
\end{itemize}

Since $Y$ avoids $[r-2]$, {we know that for a fixed $u\in V(\Gamma)\setminus Y$,} apart from at most $(r-2){4r \choose r-1}$ choices of $v \in V(\Gamma)\setminus (Y\cup\{u\})$ the edges $e_{uv}=Y\cup\{u,v\}$ of $H$ are of color $r-1$, so
$d_{\Gamma}(u)\geq n-r{4r \choose r-1}$. Also for each $1\leq i\leq r-2$, we have  $d_{\Gamma}(y_i)=n-(r-2)$. {One can easily see that Dirac's condition implies that the  graph $\Gamma$ is Hamiltonian; see \cite{Bondy}.}

Now we show that every Hamiltonian cycle in $\Gamma$ can be extended to a monochromatic Hamiltonian Berge-cycle of color $r-1$ in $H$.
Suppose that $v_{1}, v_{2},\ldots ,v_{n}$ are the vertices of a Hamiltonian cycle $C$ in $\Gamma$. {Now we define the distinct edges $f_1,f_2,\ldots,f_{n}\in E(H)$ of color $r-1$ one by one (in the same order as their subscripts appear),  such that for each $i=1,2,\ldots, n$ we have $\{v_{i},v_{i+1}\}\subseteq  f_{i}$} and $f_{1}, f_{2},\ldots ,f_{n}$  make a Hamiltonian Berge-cycle with the core sequence $v_{1}, v_{2}, \ldots ,v_{n}$.
{ We would like to choose} $f_{i}=e_{v_{i}v_{i+1}}\in F_1$ for  $v_{i}v_{i+1}\in E_{1}$. Now assume $v_{i}v_{i+1}\in E_{2}$. {Choose} $f_{i}=Y\cup\{v_{i},v_{i+1},u_i\}$ of color $r-1$ {with} $u_i\in V(\Gamma)\setminus (Y\cup\{v_{i-1},v_{i},v_{i+1},v_{i+2}\})$ and $f_i\neq f_j$ for every $j<i$. {Such an edge $f_i$ exists} since at least $n-r{4r \choose r-1}$ edges of color $r-1$ contain $Y\cup\{v_{i},v_{i+1}\}$, and $\{v_{i},v_{i+1}\}$ {lies in at most $2(r-3)+2$ edges} $f_{j}$ for $j<i$. Note that for $v_{i}v_{i+1}\in E_{2}$, if $\{v_{i},v_{i+1}\}\subseteq f_{j}$ {for some $2\leq j\leq i-2$, then $v_{j}v_{j+1}\in E_2$, $|Y\cap\{v_{i},v_{i+1},v_{j},v_{j+1}\}|=2$ and $f_j=Y\cup\{v_{i},v_{i+1},v_{j},v_{j+1}\}$. On the other hand, since each edge of $E_2$ has exactly one vertex $y_i$ for some $1\leq i\leq r-2$ we have $|E(C)\cap E_{2}|=2(r-2)$  and so $\{v_{i},v_{i+1}\}\in E_{2}$ has been used in at most $2(r-3)$ edges $f_{j}$ for $2\leq j\leq i-2$. Therefore, at most $2(r-3)+2$ edges $f_j$ for $1\leq j\leq i-1$ contain $\{v_i,v_{i+1}\}$ if $v_iv_{i+1}\in E_2$.}

\bigskip\noindent \textbf{Case 2.} $f\leq r-3.$\\
First we prove the following claim.
\begin{emp}\label{C1}
$|\overline{U}_{f+1}|\leq r-2$.
\end{emp}

Suppose to the contrary that $|\overline{U}_{f+1}|\geq r-1$. Now let $$M=\{xy_1y_2\ldots y_{f}u_{f+1}u_{f+2}\ldots u_{r-1}|u_{i}\in \overline{U}_{i}\}.$$

For each $f+1\leq i\leq r-1$, {by the definition of $\overline{U}_i$ we have $i\notin L^{*}(xy)$ for every $y\in \overline{U}_i$ and so there} are at most $(r-2)|\overline{U}_{i}|$ edges in $M$ of color $i$. {On the other hand, $\{y_i\}_{i=1}^{f}$ avoids the set of colors $\{1,2,\ldots,f\}$ and so }at most ${4r \choose r-1}$ edges in $M$ are of color $i$ for every $1\leq i\leq f$. Therefore,
{$$|M|\leq (r-2)\sum_{i=f+1}^{r-1}|\overline{U}_{i}|+f{4r \choose r-1}.$$
The inequalities (3), Remark \ref{remark1} and the assumption $|\overline{U}_{f+1}|\geq r-1$ imply that
 $$(r-1)^{r-f-3}(s-(r-f-3)(r-1))|\overline{U}_{r-1}|\leq |M|=\prod_{i=f+1}^{r-1}|\overline{U}_{i}|\leq (r-2)(s+|\overline{U}_{r-1}|)+f{4r \choose r-1},$$}
 where $s=\sum_{i=f+1}^{r-2}|\overline{U}_{i}|$. Therefore $p(s)\leq 0$, where
 $$p(x)={\big((r-1)^{r-f-3}(x-(r-f-3)(r-1))-(r-2)\big)}|\overline{U}_{r-1}|-(r-2)x-f{4r \choose r-1}.$$
Evidently, $p(x)$ is an increasing function, its derivative being positive for every $x$ { because of our assumption $|\overline{U}_{f+1}|\geq r-1$.} {By Claim \ref{claime1} we have $|\overline{U}_{r-1}|\geq (n-1)/2$. On the other hand, $s\geq (r-f-2)(r-1)$, $f\leq r-3$ and $n>6r{4r \choose r-1}$. Hence we have}
$$p(s)\geq p((r-f-2)(r-1))$$$$={\big((r-1)^{r-f-2}-(r-2)\big)}|\overline{U}_{r-1}|-(r-f-2)(r-2)(r-1)-f{4r \choose r-1}>0,$$
a contradiction. Hence $|\overline{U}_{f+1}|\leq r-2$.\\

Assume that $B_i=\overline{U}_{i}\setminus ((\cup_{j>i}\overline{U}_{j})\cup\{y_i\})$ for every $1\leq i\leq r-1$ and $\overline{U}_{f+1}=\{u_1(=y_{f+1}),u_2,\ldots,u_l\}$. According to Claim \ref{C1}, we have $l\leq r-2$. {Let $U$ be the set of all vertices $y\notin Y\cup\{x,y_{f+1}\}$, for which the edge $Y\cup\{x,y\}$ is of color $f+1$. We know that  $\{y_i\}_{i=1}^{f}$ avoids the set of colors $\{1,2,\ldots,f\}$ and so for each $1\leq i\leq f$ the number of vertices $y\notin Y\cup\{x,y_{f+1}\}$, for which the edge $Y\cup\{x,y\}$ is of color $i$ is at most ${4r \choose r-1}$. On the other hand, $i\notin L^{*}(xy_i)$ for every $f+1\leq i\leq r-1$ and so for each $f+1\leq i\leq r-1$ the number of vertices $y\notin Y\cup\{x,y_{f+1}\}$, for which the edge $Y\cup\{x,y\}$ is of color $i$ is at most $r-2$. Therefore, we have $|U|\geq n-r{4r \choose r-1}$. Let }$U$ be partitioned into $A_{1},A_{2},\ldots,A_{r-1}$, where $|A_{r-1}|=\lfloor\frac{n}{2}\rfloor+1$, $A_{f+1}=\emptyset$ and $||A_i|-|A_j||\leq 1$ for every $1\leq i,j\leq r-2$ with $i,j\neq f+1$.  Consider a graph $\Gamma$ with vertex set $V(\Gamma)=V(H)$ and edge set
$E(\Gamma)=\bigcup _{i=1}^{5}E_{i}$, where $E_i$'s are defined as follows:
\begin{itemize}
\item[$\bullet$]
$E_{1}=\{ uv\vert u\in B_i, i\neq f+1, v\notin Y\cup\{ x,u\}, c(Y_i\cup\{ x,u,v\})=f+1\}$. For each $uv\in E_{1}$, set $e_{uv}=Y_i\cup\{ x,u,v\}$, where $i$ is the minimum number such that $i\neq f+1$, $B_i\cap \{u,v\}\neq\emptyset$ and $c(Y_i\cup\{ x,u,v\})=f+1$. Now let $F_{1}=\{e_{uv}|uv\in E_{1}\}$.

Note that for every $1\leq i\leq r-1$ we have $B_i=\overline{U}_{i}\setminus ((\cup_{j>i}\overline{U}_{j})\cup\{y_i\})$ and by Fact \ref{ui}, $B_i\cap Y=\emptyset$. Therefore in the subgraph of $G=\Gamma(H)$ induced by the edges $E_1$ the vertices $Y$ are isolated vertices. Now we define the edges crossing the vertices $Y$.

\item[$\bullet$]
$E_{2}=\{ y_iv\vert v\in A_i, 1\leq i\leq r-1, i\neq f+1\}$ and for each $y_iv\in E_{2}$, set $e_{y_iv}=Y\cup\{x,v\}$. Also, let $F_{2}=\{e_{y_iv}|y_iv\in E_{2}\}$.
\item[$\bullet$]
{Now we define new edges to increase the degrees of vertices in $\overline{U}_{f+1}$ with small degrees in the subgraph of $G=\Gamma(H)$ induced by the edges $E_1\cup E_2$. In fact we define a set of new edges $E_3$ such that the degree of each vertex $u_i$ for $1\leq i\leq l$ in the subgraph of $G=\Gamma(H)$ with vertex set $V(H)$ and edge set $E_1\cup E_2\cup E_3$ is at least $2r+1$. For each $1\leq i\leq l$ we have $d_{f+1}(u_i)>{4r \choose r-1}$, since otherwise $\{u_i, y_1,y_2,\ldots,y_f\}$ avoids the set of colors $[f+1]$, which is a contradiction to the maximality of $f$. We use this fact here.} To define $E_3$, we do the following. Let $\Gamma_1$ be the graph with vertex set $V(H)$ and edge set $E_{1}\cup E_{2}$. For each $1\leq i\leq l$ assume that $\overline{N}_i=Y\cup\overline{U}_{f+1}\cup N_{\Gamma_1}(u_i)\cup\{x\}$ and set $t_i=0$ if $d_{\Gamma_1}(u_{i})>2r$ and $t_i=2r+1-d_{\Gamma_1}(u_{i})$, otherwise. Now we show that there are $\sum_{i=1}^{l}t_i$ distinct edges $e_{ij}\notin F_{1}\cup F_{2}$ (where $1\leq i\leq l$ and $1\leq j\leq t_i$) of color $f+1$ with $u_i\in e_{ij}$ such that for each $1\leq i\leq l$ there exist $t_i$ distinct vertices $v_{ij}\in e_{ij}\setminus \overline{N}_i$. For this, set $r_{11}=0$, $N_{11}=\overline{N}_1$ and $E_{11}=F_{1}\cup F_{2}$ and follow the following step for $i=1,2,\ldots,l$ if $t_i>0$.

{\bf Step i:} For each $1\leq j\leq t_i$, since $d_{f+1}(u_i)>{4r \choose r-1}\geq {|N_{ij}|-1\choose r-1}+r_{ij}$, there is an edge $e_{ij}\notin E_{ij}$ of color $f+1$ which contains $u_i$ and a vertex $v_{ij}\in e_{ij}\setminus N_{ij}$. Note that since $\{y_i\}_{i=1}^{f}$ avoids $[f]$ and $f$ is maximum subject to this property, we have $d_{f+1}(u_i)>{4r \choose r-1}$. Now set $r_{i(j+1)}=r_{ij}+1$, $N_{i(j+1)}=N_{ij}\cup\{v_{ij}\}$ and $E_{i(j+1)}=E_{ij}\cup \{e_{ij}\}$ and continue the above procedure. We apply the above procedure $t_i$ times to find the edges $e_{ij}$ and the vertices $v_{ij}$ for $1\leq j\leq t_i$ with desired properties. Finally let $r_{(i+1)1}=r_{i(t_i+1)}$, $N_{(i+1)1}=\overline{N}_{i+1}$ and $E_{(i+1)1}=E_{i(t_i+1)}$ and go to Step $i+1$.

Clearly, $E_{l(t_l+1)}\setminus E_{11}$ contains $\sum_{i=1}^{l}t_i$ distinct edges $e_{ij}$ with desired properties. Now set $A=\bigcup_{i=1}^{l}\bigcup_{j=1}^{t_i}e_{ij}$, $\overline{E}_i=\{u_{i}v_{ij}|1\leq j\leq t_i\}$, $\overline{F}_i=\{e_{ij}|1\leq j\leq t_i\}$, $E_{3}=\bigcup_{i=1}^{l}\overline{E}_i$ and $F_{3}=\bigcup_{i=1}^{l}\overline{F}_i$.

\item[$\bullet$]
{The set of edges $E_4$ is defined in a more or less similar way. Here we define these edges to increase the degrees of vertices in $U_{\{1,2,\ldots,r-1\}}$ with small degrees in the subgraph of $G=\Gamma(H)$ induced by the edges $E_1\cup E_2\cup E_3$, where $U_{\{1,2,\ldots,r-1\}}=\bigcap_{i=1}^{r-1}U_i$.  In fact we define a set of new edges $E_4$ such that the degree of each vertex in $U_{\{1,2,\ldots,r-1\}}$ in the subgraph of $G$ with vertex set $V(H)$ and edge set $\bigcup_{i=1}^{4}E_i$ is at least $2r+1$. We will see this result in Fact \ref{Factd4}. To define $E_4$, we do the following:} Assume that $U_{\{1,2,\ldots,r-1\}}=\{w_{1}, w_{2},\ldots, w_{m}\}$ and $d_{\Gamma_2}(w_{1})\leq d_{\Gamma_2}(w_{2})\leq \cdots \leq d_{\Gamma_2}(w_{m})$, where $\Gamma_2$ is the graph with vertex set $V(H)$ and edge set $\bigcup _{i=1}^{3} E_{i}$.
For each $1\leq i\leq r'=\min\{r,m\}$, set $t'_i=0$ when $d_{\Gamma_2}(w_{i})>2r$. Otherwise, set $t'_i=2r+1-d_{\Gamma_2}(w_{i})$. Also, set $N'_i=Y\cup\overline{U}_{f+1}\cup N_{\Gamma_2}(w_i)\cup\{x\}$.  An argument similar to the one used in the definition of $E_3$ shows that there are $\sum_{i=1}^{r'}t'_i$  distinct edges $e'_{ij}\notin F_{1}\cup F_{2}\cup F_{3}$ (where $1\leq i\leq r'$ and $1\leq j\leq t'_i$) of color $f+1$ with $w_i\in e'_{ij}$ such that for each $1\leq i\leq r'$ there exist $t'_i$ distinct vertices $v'_{ij}\in e'_{ij}\setminus N'_i$. Now, set $B=\bigcup_{i=1}^{r'}\bigcup_{j=1}^{t'_i}e'_{ij}$, $E'_i=\{w_{i}v'_{ij}|1\leq j\leq t'_i\}$, $F'_i=\{e'_{ij}|1\leq j\leq t'_i\}$, $E_{4}=\bigcup_{i=1}^{r'}E'_i$ and $F_{4}=\bigcup_{i=1}^{r'}F'_i$.
\item[$\bullet$]
$E_{5}=\{ xv\vert v\in V(\Gamma)\setminus (Y\cup \overline{U}_{f+1}\cup A\cup B)\}$.
\end{itemize}

{In the following fact using the above definitions we see that the set of edges $F_1,F_1,F_1,F_1$ are pairwise disjoint.}

{\begin{fact}\label{disjoint-Fi}
For each $1\leq i,j\leq 4$ and $i\neq j$, we have $F_i\cap F_j=\emptyset$.
\end{fact}}

{First we show that $F_1\cap F_2=\emptyset$. With a contrary assume that $f\in F_1\cap F_2$. Since $f\in F_1$ from the definition of $F_1$ we have $f=e_{uv}=Y_i\cup\{ x,u,v\}$, where $u\in B_i$, $i\neq f+1$, $v\notin Y\cup\{ x,u\}$ and $c(Y_i\cup\{ x,u,v\})=f+1$. One can easily see that $y_i\notin f$. On the other hand, $f\in F_2$. Hence $f=e_{y_jz}=Y\cup\{x,z\}$ for some $z\in A_j$, where $1\leq j\leq r-1$ and $j\neq f+1$. Hence $y_i\in Y\subseteq f$, a contradiction.  Therefore, $F_1\cap F_2=\emptyset$. Now using the definition of $E_3$, we have $F_{3}=\bigcup_{i=1}^{l}\overline{F}_i$ and $\overline{F}_i=\{e_{ij}|1\leq j\leq t_i\}$. On the other hand, $e_{ij}\notin F_{1}\cup F_{2}$ for every $1\leq i\leq l$ and $1\leq j\leq t_i$. Therefore, $F_3\cap (F_{1}\cup F_{2})=\emptyset$.
Again from the definition of $E_4$, we have $F_{4}=\bigcup_{i=1}^{r'}F'_i$ and $F'_i=\{e'_{ij}|1\leq j\leq t'_i\}$. Moreover, $e'_{ij}\notin F_{1}\cup F_{2}\cup F_{3}$ for every $1\leq i\leq r'$ and $1\leq j\leq t'_i$. Therefore, $F_4\cap (F_{1}\cup F_{2}\cup F_{3})=\emptyset$.
}

\begin{emp}\label{Hhamilt2}
The graph  $\Gamma$ is Hamiltonian.
\end{emp}

{\it Proof of Claim \ref{Hhamilt2}:}
Assume that $d_1\leq d_2\leq \cdots\leq d_n$ are degrees of the vertices of $\Gamma$.  We show that $d_1>2r$ and $d_{n-i}\geq n-i$ for each $2r-1 \leq i\leq \frac{n}{2}$. Therefore, Lemma \ref{Ch1} implies the existence of a Hamiltonian cycle in $\Gamma$. Now we give the following facts about the degrees of vertices of $\Gamma$.

\begin{fact}\label{Factd1}
$d_{\Gamma}(x)\geq n-4r^3$.
\end{fact}

{To see Fact \ref{Factd1} note that using the definitions  $A$ and $B$ (in the definitions of $E_3$ and $E_4$) and Claim \ref{C1} (that indicates $l\leq r-2$) and the fact $r'\leq r$, we have $$|A|\leq r(t_1+t_2+\cdots+t_l)\leq r(2r+1)l\leq r(r-2)(2r+1),$$
and
$$|B|\leq r(t'_1+t'_2+\cdots+t'_{r'})\leq r^2(2r+1).$$
Therefore,  $d_{\Gamma}(x)=n-|Y\cup \overline{U}_{f+1}\cup A\cup B|\geq n-4r^3$.}

\begin{fact}\label{Factd2}
For each $1\leq i\leq r-1$ with $i\neq f+1$ and each $u\in \overline{U}_{i}\setminus \{y_{i}\}$ we have
 $d_{\Gamma}(u)>n-r{4r \choose r-1}$. Moreover, for every $u\in \overline{U}_{f+1}$, we have $d_{\Gamma}(u)>2r.$
\end{fact}

To show Fact \ref{Factd2} note that Fact \ref{avoid1} implies that the set of vertices $Y_i\cup\{ x\}$ avoids the set of colors $[r-1]\setminus\{i,f+1\}$. On the other hand, $i\notin L^{*}(xu)$ for $u\in \overline{U}_{i}\setminus \{y_{i}\}$ and so  $(Y_i\cup\{ x,u\})$ avoids all colors $[r-1]\setminus\{f+1\}$. Therefore, apart from at most $(r-2){4r \choose r-1}$ choices of $v \in V(\Gamma)\setminus (Y\cup\{ x,u\})$ we have $uv\in E_1$ and so $d_{\Gamma}(u)>n-r{4r \choose r-1}$. Moreover, for every $u_i\in \overline{U}_{f+1}$, we have $d_{\Gamma}(u_{i})\geq d_{\Gamma_1}(u_{i})+t_i>2r$ (see the definition $E_3$).\\

\begin{fact}\label{Factd3}
$d_{\Gamma}(y_{r-1})>n/2$ and $d_{\Gamma}(y_{i})>2r$ for each $1\leq i\leq r-1$ and $i\neq f+1$.
\end{fact}
Fact \ref{Factd3} follows from the fact  $y_{i}v\in E(\Gamma)$ for each $v\in A_{i}$ and $|A_{r-1}|>n/2$ and $|A_{i}|>2r$ for each $1\leq i\leq r-1$ and $i\neq f+1$.

\begin{fact}\label{Factd4}
 $d_{\Gamma}(u)>2r$ for each $u\in U_{\{1,2,\ldots,(r-1)\}}$.
\end{fact}

To see Fact \ref{Factd4} assume  $U_{12\ldots (r-1)}=\{w_{1}, w_{2},\ldots, w_{m}\}\neq \emptyset$. We claim that $$\min\{d_{\Gamma}(w_{i})| 1\leq i\leq m\}>2r.$$
First assume that $m\leq r$. According to the definition of $E_4$, for each $1\leq i\leq m$, we have $d_{\Gamma}(w_{i})\geq d_{\Gamma_2}(w_{i})+t'_i>2r$, where $\Gamma_2$ is the graph  with vertex set $V(\Gamma)$ and edge set $\bigcup _{i=1}^{3} E_{i}$.
Now let $m\geq r+1$, $\vert \overline{U}_{r-1}\setminus\{y_{r-1}\}\ \vert =k$ and $d_{\Gamma_2}(w_{1})\leq d_{\Gamma_2}(w_{2})\leq \cdots \leq d_{\Gamma_2}(w_{m})$. Again, according to the definition of the edges $E_4$, we have $d_{\Gamma}(w_i)> 2r$ for $1\leq i\leq r$ and so it suffices to show that  $d_{\Gamma}(w_{r+1})\geq d_{\Gamma_2}(w_{r+1})>2r$. For $i=1,\ldots ,m$, consider
$$N_{i}=\{\{ x,y_1,y_2,\ldots,y_{r-2},v,w_{i}\}\setminus\{y_{f+1}\}\ \vert v\in \overline{U}_{r-1}\setminus\{y_{r-1}\}\}.$$
For every $1\leq i\leq m$, suppose that $n_{i}$ is the number of edges of color $f+1$ in $N_{i}$.
Clearly for each $1\leq i\leq m$, the edges of color $f+1$ in $N_i$ belong to $F_1$ and so $d_{\Gamma_2}(w_{i}) \geq n_{i}$. Moreover, the vertices $\{ x,y_1,y_2,\ldots,y_{r-2}\}\setminus\{y_{f+1}\}$ avoids the colors $[r-1]\setminus \{f+1,r-1\}$ and $r-1\notin L^{*}(xv)$ for each $v\in \overline{U}_{r-1}\setminus\{y_{r-1}\}$. Therefore, among all $mk$ edges in $\bigcup _{i=1}^{m}N_{i}$, there are at most ${4r \choose r-1}$ edges of color $i$ for each $i\neq f+1, r-1$ and $(r-2)k$ edges of color $r-1$. Thus, $$\sum_{i=1}^{m}n_{i}\geq (m-r+2)k-(r-3){4r \choose r-1}.$$
If $d_{\Gamma_2}(w_{r+1})\leq 2r$,  then $$\sum_{i=1}^{r+1} n_{i}\leq \sum_{i=1}^{r+1} d_{\Gamma_2}(w_{i})\leq 2r(r+1).$$ Therefore
$$\sum_{i=r+2}^{m}n_{i} \geq (m-r+2)k-(r-3){4r \choose r-1}-2r(r+1)>(m-r+1)k,$$
which is impossible since $\vert \bigcup _{i=r+2}^{m}N_{i}\vert=(m-r-1)k$.
Thus  $d_{\Gamma}(w_{r+1})\geq d_{\Gamma_2}(w_{r+1})>2r$ and consequently $d_{\Gamma}(w_{i})>2r$  for  $r+1\leq i\leq m$.
On the other hand, according to the definition of $\Gamma$, we have $d_{\Gamma}(w_{i})\geq d_{\Gamma_2}(w_{i})+t'_i>2r$ for each  $1\leq i\leq r$ and so $\min\{d_{\Gamma}(w_{i})| 1\leq i\leq m\}>2r$.\\

Clearly $V(H)=V(\Gamma)=(\cup_{i=1}^{r-1} \overline{U}_{i})\cup \{y_i\}_{i=1}^{f}\cup U_{\{1,2,\ldots,(r-1)\}}\cup \{x\}.$ Therefore, Facts \ref{Factd1}-\ref{Factd4} imply that the minimum degree of $\Gamma$ is greater than $2r$ and so $d_1>2r$. Now we are going to show that $d_{n-i}\geq n-i$ for each $2r-1 \leq i\leq \frac{n}{2}$. To see this, first we show that most of the vertices of $\overline{U}_{r-1}$ have degree greater than $n-2r$ in $\Gamma$. For this, let $D_i$ be the set of all edges of color $i$ containing the vertices of $Y_{r-1}\cup \{x\}$ for each $i\neq f+1, r-1$ and let $$W=\bigcup_{i\neq f+1, r-1}\bigcup_{e\in D_i}(e\setminus (Y_{r-1}\cup \{x\})).$$ Using Fact \ref{avoid1}, $Y_{r-1}\cup\{ x\}$ avoids each color $i\neq f+1, r-1$, hence $|D_i|\leq {4r \choose r-1}$. On the other hand, for each $i\neq f+1,r-1$ and each $e\in D_i$ we have $|e\setminus (Y_{r-1}\cup \{x\})|=2$ and so $|W|\leq 2(r-3){4r \choose r-1}$. For every $u\in \overline{U}_{r-1}\setminus (W\cup \{y_{r-1}\})$, $r-1\notin L^{*}(xu)$ and so we have $uv\in E_1$, apart from at most $r-2$ choices of $v \in V(\Gamma)\setminus (Y\cup\{ x,u\})$. Moreover,
 for every $u\in \overline{U}_{r-1}\cap W\setminus \{y_{r-1}\}$, apart from at most $(r-2){4r \choose r-1}$ choices of $v \in V(\Gamma)\setminus (Y\cup\{ x,u\})$ we have $uv\in E_1$ and so
 $d_{\Gamma}(u)>n-r{4r \choose r-1}$. Hence we have the following fact.
\begin{fact}\label{Factd5}
$d_{\Gamma}(u)>n-2r$, where $u\in \overline{U}_{r-1}\setminus (W\cup \{y_{r-1}\})$. Moreover, for each $u\in \overline{U}_{r-1}\cap W\setminus \{y_{r-1}\}$, we have $d_{\Gamma}(u)>n-r{4r \choose r-1}$.
\end{fact}

By Fact \ref{Factd5} for each vertex $u\in \overline{U}_{r-1}\setminus (W\cup \{y_{r-1}\})$, we have $d_{\Gamma}(u)>n-2r$. Moreover, since $|\overline{U}_{r-1}|\geq (n-1)/2$ and $|W|\leq 2(r-3){4r \choose r-1}$ we have $$|\overline{U}_{r-1}\setminus (W\cup \{y_{r-1}\})|\geq\frac{n-3}{2}-2(r-3){4r \choose r-1},$$ and so at least $\lceil\frac{n-3}{2}\rceil-2(r-3){4r \choose r-1}$ vertices of $\Gamma$ has degree greater than $n-2r$, this means that

$$~~~~~~~~~~~~~~~~~~~~~~~~~~~~~~~i\geq \lfloor \frac{n+5}{2}\rfloor+2(r-3){4r \choose r-1}\Longrightarrow d_{i}>n-2r~~~~~~~~~~~~~~~~~~~~~~~~~~~~~(4)$$

 Fact \ref{Factd2} implies that for each $1\leq i\leq r-1$ and $i\neq f+1$  and for every  $u\in \overline{U}_{i}\setminus \{y_{i}\}$,  we have
 $d_{\Gamma}(u)>n-r{4r \choose r-1}$. Now using Facts \ref{Factd1}, we have $d_{\Gamma}(x)\geq n-4r^3$. On the other hand, $|\overline{U}_{r-1}|\geq (n-1)/2$ and $n-4r^3>n-r{4r \choose r-1}$ and so at least $\lceil \frac{n-1}{2}\rceil$ vertices of $\Gamma$ have degree greater than $n-r{4r \choose r-1}$, this means that

 $$~~~~~~~~~~~~~~~~~~~~~~~~~~~~~~~i\geq\lfloor\frac{n+3}{2}\rfloor\Longrightarrow d_{i}>n-r{4r \choose r-1}~~~~~~~~~~~~~~~~~~~~~~~~~~~~~(5)$$

 Now using Fact \ref{Factd3}, we have $d_{\Gamma}(y_{r-1})>n/2$. Therefore, we have

 $$~~~~~~~~~~~~~~~~~~~~~~~~~~~~~~~i\geq\lfloor  \frac{n+1}{2}\rfloor\Longrightarrow d_{i}>n/2~~~~~~~~~~~~~~~~~~~~~~~~~~~~~(6)$$

Since $n>6r{4r \choose r-1}$ using (4),(5),(6) we conclude that $d_{n-i}\geq n-i$ for each $2r-1 \leq i\leq \frac{n}{2}$. On the other hand, $d_1>2r$. Now, Lemma \ref{Ch1} implies the existence of a Hamiltonian cycle in $\Gamma$.
\epf
\begin{emp}\label{extend2}
There is a monochromatic Hamiltonian Berge-cycle of color $f+1$ in $H$.
\end{emp}
{\it Proof of Claim \ref{extend2}:}
We show that every Hamiltonian cycle in $\Gamma$ can be extended to a monochromatic Hamiltonian Berge-cycle of color $f+1$ in $H$.
Suppose that $v_{1},v_{2},\ldots,v_{n}=x$ are the vertices of a Hamiltonian cycle $C$ in $\Gamma$. {Now for $i=1,2,\ldots, n$, we define the edges $f_{i}\in E(H)$ of color $f+1$ one by one (in the same order as their subscripts appear), so that $\{v_{i},v_{i+1}\}\subseteq  f_{i}$ and $f_{1}, f_{2},\ldots, f_{n}$ make a Hamiltonian Berge-cycle with the core sequence $v_{1}, v_{2}, \ldots ,v_{n}$.
First we follow the following step for $i=1,2,\ldots, n-2$ one by one to define the edges $f_{1}, f_{2},\ldots, f_{n-2}$.}\\

{{\bf Step i:} If $v_{i}v_{i+1}\in E_{j}$ for some  $j\in\{1,2\}$, then set $f_{i}=e_{v_{i}v_{i+1}}\in F_{j}$. Let $f_{i}=e_{kj}\in F_{3}$ if $\{v_{i},v_{i+1}\}=\{u_k,v_{kj}\}$ and $u_kv_{kj}\in E_{3}$, where $k\in\{1,2,\ldots,l\}$ and $1\leq j\leq t_k$. Finally, let $f_{i}=e'_{kj}\in F_{4}$ if $\{v_{i},v_{i+1}\}=\{w_k,v'_{kj}\}$ and $w_kv'_{kj}\in E_{4}$, where $k\in\{1,2,\ldots,r'\}$ and $1\leq j\leq t'_k$. Then go to Step $i+1$.}\\

{According to the definitions of $F_1,F_2,F_3$ and $F_4$, for each $1\leq i\leq n-2$ the edge $f_{i}\in \bigcup_{i=1}^{4} F_i$ is of color $f+1$ and $\{v_{i},v_{i+1}\}\subseteq  f_{i}$. Now we claim that $f_i\neq f_j$ for every $i\neq j$ with $1\leq i,j\leq n-2$. It only suffices to prove the following fact.}

{\begin{fact}\label{Final}
For each $1\leq i\leq n-2$ and $1\leq j< i$, we have $f_i\neq f_j$.
\end{fact}}

{Assume that $f_i\in F_{r_i}$ and $f_j\in F_{r_j}$, where $r_i,r_j\in\{1,2,3,4\}$. Using Fact \ref{disjoint-Fi}, $F_{r_i}\cap F_{r_j}=\emptyset$ if $r_i\neq r_j$. Hence $f_i\neq f_j$ when $r_i\neq r_j$. Therefore, we may assume that $r_i=r_j$. First assume that $j=i-1$. We divide our proof for this fact into some cases:\\}

{First let $r_{i-1}=r_i=1$. Then $f_{i-1}=e_{v_{i-1}v_{i}}=Y_p\cup\{ x,v_{i-1},v_{i}\}$ and $f_{i}=e_{v_{i}v_{i+1}}=Y_q\cup\{ x,v_{i},v_{i+1}\}$, where $p,q$ are the minimum numbers such that $p,q\neq f+1$, $B_p\cap \{v_{i-1},v_{i}\}\neq\emptyset$, $B_q\cap \{v_{i},v_{i+1}\}\neq\emptyset$ and $c(Y_p\cup\{ x,v_{i-1},v_{i}\})=c(Y_q\cup\{ x,v_{i},v_{i+1}\})=f+1$. One can easily see that $\{v_{i-1},v_{i}\}\nsubseteq f_i$ and so $f_i\neq f_{i-1}$.}

{Now let $r_{i-1}=r_i=2$. Then $\{v_{i-1},v_{i}\}=\{y_t,v\}$ for some $1\leq t\leq r-1$, $t\neq f+1$, $v\in A_t$ and $f_{i-1}=e_{v_{i-1}v_{i}}=e_{y_tv}=Y\cup\{x,v\}$. Since $A_p\cap A_q=\emptyset$ for $p\neq q$ and $r_i=2$, we have $v_i=y_t$, $v_{i-1}, v_{i+1}\in A_t$ and $f_i=e_{v_{i}v_{i+1}}=e_{y_tv_{i+1}}=Y\cup\{x,v_{i+1}\}$. Here clearly $v_{i+1}\notin f_{i-1}$ and so $f_i\neq f_{i-1}$.}

{Now let $r_{i-1}=r_i=3$, then by the definitions of $E_3$ and $F_3$ in Page 9 we have $f_{i-1}=e_{k_1j_1}\in F_{3}$ and $f_{i}=e_{k_2j_2}\in F_{3}$, where $\{v_{i-1},v_{i}\}=\{u_{k_1},v_{k_1j_1}\}$ and $\{v_{i},v_{i+1}\}=\{u_{k_2},v_{k_2j_2}\}$ for some $k_1, k_2\in\{1,2,\ldots,l\}$, $1\leq j_1\leq t_{k_1}$ and $1\leq j_2\leq t_{k_2}$. Now assume with a contrary that $f_{i-1}=f_i$. Using the definitions of $E_3$ and $F_3$, we have $v_{k_1j_1}, v_{k_2j_2}\notin \overline{U}_{f+1}$ and so $v_i=u_{k_1}=u_{k_2}$, $k_1=k_2$, $v_{i-1}=v_{k_1j_1}$ and $v_{i+1}=v_{k_2j_2}$. On the other hand, $v_{i-1}\neq v_{i+1}$ and so $j_1\neq j_2$. Hence from the definition of $F_3$, we have $e_{k_1j_1}\neq e_{k_1j_2}$ and so $f_{i-1}\neq f_i$, a contradiction to our assumption.}

{Finally let $r_{i-1}=r_i=4$, then using the definitions of $E_4$ and $F_4$ in Page 10 we have $f_{i-1}=e'_{k_1j_1}\in F_{4}$ and $f_{i}=e'_{k_2j_2}\in F_{4}$, where $\{v_{i-1},v_{i}\}=\{w_{k_1},v'_{k_1j_1}\}$ and $\{v_{i},v_{i+1}\}=\{w_{k_2},v'_{k_2j_2}\}$ for some $k_1, k_2\in\{1,2,\ldots,r'\}$, $1\leq j_1\leq t'_{k_1}$ and $1\leq j_2\leq t'_{k_2}$. With the same argument we can see that $k_1\neq k_2$ or $j_1\neq j_2$. Therefore, from the definition of $F_4$, we have $e'_{k_1j_1}\neq e'_{k_1j_2}$ and so $f_{i-1}\neq f_i$.\\}

{Now assume $j\leq i-2$. In this case by the definitions of $F_1, F_2, F_3$ and $F_4$, one can easily see that $\{v_{i},v_{i+1}\}\nsubseteq f_j$ or $\{v_{j},v_{j+1}\}\nsubseteq f_i$ and so again $f_{i}\neq f_j$.\\}

{Now we are going to give the definitions of $f_{n-1}$ and $f_n$ with desired properties. First let $i=n-1 $. Since $\{ v_{n-1},x\}$  has been used in at most one of the edges  $f_{i}$'s, where $1\leq i\leq n-2$  (only possibly in $f_{n-2}$) and  $f+1\in L^{*}(v_{n-1}x)$, then  we can choose an appropriate edge $f_{n-1}$ of color $f+1$, where $f_{n-1}\neq f_i$ for each $1\leq i\leq n-2$. Similarly, for $i=n$, since  $\{x, v_{1}\}$  has been used in at most two edges  $f_{i}$'s, where $1\leq i\leq n-1$ (only possibly in $f_{1}$ and $f_{n-1}$) and  $f+1\in L^{*}(xv_{1})$, then  we can choose an appropriate edge $f_{n}$ of color $f+1$, where $f_{n}\neq f_i$ for each $1\leq i\leq n-1$.}
\epf
This finishes the proof of Theorem \ref{main result}.
\epf

\end{document}